\documentclass{mathscan}
\usepackage{amssymb}
\usepackage{enumerate}
\usepackage{graphicx}
\usepackage{amscd}
\usepackage{enumerate}
\usepackage{amssymb}
\usepackage{amsthm}
\usepackage{amsfonts}
\usepackage{fancyhdr}
\usepackage{mathrsfs}
\usepackage{amsxtra}
\usepackage{setspace}

\usepackage{color}

\newtheorem{thm}{Theorem}[section]
  
  \newtheorem{cor}[thm]{Corollary}
  \newtheorem{lem}[thm]{Lemma}

  \theoremstyle{definition}
  \newtheorem{rem}[thm]{Remark}

  \numberwithin{equation}{section}

\def\beq{\begin{eqnarray*}}
\def\eeq{\end{eqnarray*}}

\def\DD{{\mathbb{D}}}

\def\HH{{\mathcal H}}

\def\LL{{\mathcal{L}}}
\def\MM{{\mathcal M}}

\def\ker{\text{\rm ker}}

\newcommand{\mult}{\mathfrak{M}\hbox{\rm{ult}}}

\def\T{{\mathbb T}}
\def\H{{\mathscr H}}

\begin{document}

\title[Backward shift invariant subspaces in RKHS]{Backward shift invariant subspaces in reproducing kernel Hilbert spaces}

\author[Fricain]{Emmanuel Fricain}
\address{Laboratoire Paul Painlev\'e, Universit\'e Lille 1, 59 655 Villeneuve d'Ascq C\'edex }
\email{emmanuel.fricain@math.univ-lille1.fr}

\author[Mashreghi]{Javad Mashreghi}
\address{D\'epartement de math\'ematiques et de statistique, Universit\'e Laval, Qu\'ebec, QC,
Canada G1K 7P4}
\email{Javad.Mashreghi@mat.ulaval.ca}

\author[Rupam]{Rishika Rupam}
\address{Laboratoire Paul Painlev\'e, Universit\'e Lille 1, 59 655 Villeneuve d'Ascq C\'edex}
\email{Rishika.Rupam@math.univ-lille1.fr}

\thanks{The first and third authors were supported by Labex CEMPI (ANR-11-LABX-0007-01) and the grant ANR-17-CE40
-0021 of the French National Research Agency ANR (project Front). The second author was supported by grants from NSERC (Canada).}



\begin{abstract}
In this note, we describe the backward shift invariant subspaces for an abstract class of reproducing kernel Hilbert spaces. Our main result is inspired by a result of Sarason concerning de Branges--Rovnyak spaces (the non-extreme case). Furthermore, we give new applications in the context of the range space of co-analytic Toeplitz operators and sub-Bergman spaces.
\end{abstract}

\maketitle

\section{Introduction}
A celebrated theorem of Beurling describes all (non-trivial) closed invariant subspaces of the Hardy space $H^2$ on the open unit disc $\mathbb D$ which are invariant with respect to the backward shift operator $S^*$. They are of the form $K_\Theta=(\Theta H^2)^\perp$, where $\Theta$ is an inner function. The result of Beurling was the cornerstone of a whole new direction of research lying at the interaction between operator theory and complex analysis. It was generalized in many ways. See for instance \cite{MR566739,MR2244001,MR1903737,MR0270196,MR2034817}. Sarason \cite{Sarason-86-OT} classified the non-trivial closed backward shift invariant subspaces of the de Branges--Rovnyak spaces $\H(b)$, where $b$ is a non-extreme point of the closed unit ball of $H^\infty$: they are of the form $K_\Theta\cap\H(b)$, where $\Theta$ is an inner function. In other words, the closed invariant subspaces for $S^*|\H(b)$ are the trace on $\H(b)$ of the closed invariant subspaces for $S^*$. This naturally leads to the following question: \

(Q): \emph{let $\H_1$ and  $\H_2$ be two reproducing kernel Hilbert spaces on $\mathbb D$ such that $\H_1\subset\H_2$; assume that the shift operator $S$ (multiplication by the independent variable) is contractive on $\H_2$ and if $T=S|\H_2$, its adjoint $T^*$ maps $\H_1$ contractively into itself. Then, is it true that every closed invariant subspace $\mathcal E$ of $T^*|\H_1$ has the form $E\cap \H_1$, where $E$ is a closed invariant subspace of $T^*$ (as an operator on $\H_2$)? In other words, are the closed invariant subspaces for $T^*|\H_1$ the trace on $\H_1$ of the closed invariant subspaces for $T^*$?}\

It should be noted that, of course, the interesting situation is when $\H_1$ is not a closed subspace of $\H_2$. Sarason's result says that the answer to question (Q) is affirmative in the situation where $\H_2=H^2$ and $\H_1=\H(b)$, with $b$ a non-extreme point of the closed unit ball of $H^\infty$. However, it should be noted that question (Q) has a negative answer in the case where $\H_1=\mathcal D$ is the Dirichlet space and $\H_2=H^2$ is the Hardy space. Indeed, let $(z_n)_{n\geq 1}$ be a non-Blaschke sequence of $\mathbb D$ which is a zero set for $A^2$ and put
$$
M:=\{f\in A^2:f(z_n)=0, n\geq 1\},
$$
where $A^2$ is the Bergman space of $\mathbb D$. Define 
$$
N=\{F\in\mbox{Hol}(\mathbb D): F'=f,\,f\in A^2\ominus M\}.
$$
It is not difficult to see that $N$ is a non-trivial closed subspace of $\mathcal D$, which is $S^*$ invariant. Then, observe that $N$ cannot be of the form $E\cap \mathcal D$, where $E$ is a closed subspace of $H^2$ invariant with respect to $S^*$. Indeed, assume on the contrary that there exists a closed subspace $E$ of $H^2$, invariant with respect to $S^*$, such that $N=E\cap \mathcal D$. Since $N$ is non-trivial, the subspace $E$ is also non-trivial, and by Beurling's theorem, there exists an inner function $\Theta$ such that $E=K_\Theta$. Thus $N=K_\Theta\cap \mathcal D$. Observe now that for every $n\geq 1$, the Cauchy kernel $k_{\lambda_n}$ belongs to $N$ (because its derivative is up to a constant the reproducing kernel of $A^2$ at point $\lambda_n$ and thus it is orthogonal to $M$). Then $k_{\lambda_n}\in K_\Theta$, $n\geq 1$. To get a contradiction, it remains to see that, since $(z_n)_{n\geq 1}$ is not a  Blaschke sequence, then the sequence of Cauchy kernels $k_{\lambda_n}$, $n\geq 1$, generates all $H^2$, and we deduce that $H^2\subset K_\Theta$, which is absurd. 

The aim of this note is to present a general framework where the answer to the question (Q) is affirmative. Note that in \cite{Aleman-Malman}, Aleman--Malman present another general situation of reproducing kernel Hilbert spaces where they extend Sarason's result. 

In Section 2, we first recall some basic facts on reproducing kernel Hilbert spaces and on the Sz.-Nagy--Foias model for contractions. Then, in Section 3, we study the properties of multiplication operators in our general context and prove that the scalar spectral measures of the minimal unitary dilation of $T^*|\H_1$ are absolutely continuous. We also show that when $\mathcal H_2=H^2$, then the reproducing kernel Hilbert space $\H_1$ satisfies an interesting division property, the so-called $F$-property. In Section 4, we give an analogue of Beurling's theorem in our general context and give an application to cyclic vectors for the backward shift. 
In Section 5, we show that our main theorem can be applied to $\H(b)$ spaces and range space of co-analytic Toeplitz operators. We also provide a new application in the context of sub-Bergman Hilbert space which was recently studied in \cite{sultanic2006sub,zhu1996sub,MR1990528}.

\section{Preliminaries}
We first recall some standard facts on reproducing kernel Hilbert spaces. See \cite{Paulsen} for a detailed treatment of RKHS.
\subsection{Reproducing kernel Hilbert spaces and multipliers}
Let $\H$ be a Hilbert space of complex valued functions on a set $\Omega$. We say that $\H$ is a {\em reproducing kernel Hilbert space} (RKHS) on $\Omega$ if the following two conditions are satisfied:
\begin{enumerate}
\item[(P1)] for every $\lambda\in \Omega$, the point evaluations $f \longmapsto f(\lambda)$ are bounded on $\H$;
\item[(P2)] for every $\lambda\in \Omega$, there exists a function $f\in\H$ such that $f(\lambda)\neq 0$.
\end{enumerate}
According to the Riesz representation theorem, for each $\lambda \in \Omega$, there is a function $k^{\H}_\lambda$ in $\H$, called the  {\em reproducing kernel} at point $\lambda$, such that
\[
f(\lambda) = \langle f,k^{\H}_\lambda \rangle_{\H}, \qquad (f \in \H).
\]
Note that according to (P2), we must have $k_\lambda^\H\not\equiv 0$. Moreover if $(f_n)_n$ is a sequence in $\H$, then
\begin{equation}\label{eq:weakly-convergence-simply-convergence}
f_n\to f \mbox{ weakly in }\H\implies \forall \lambda\in \Omega,\,\, \lim_{n\to \infty}f_n(\lambda)=f(\lambda).
\end{equation}

A \emph{multiplier} of $\H$ is a complex valued function $\varphi$ on $\Omega$ such that $\varphi f \in \H$ for all $f \in \H$. The set of all multipliers of $\H$ is denoted by $\mult(\H)$. Using the closed graph theorem, we see that if $\varphi$ belongs to $\mult(\H)$, then the map
\begin{equation}\label{eq:multiplier-defn}
M_{\varphi,\H}:\left |\begin{array}{cccl}
&\H& \longrightarrow &\H\cr
& f &\longmapsto &\varphi f 
\end{array} \right.
\end{equation}
is  bounded on $\H$. When there is no ambiguity, we simply write $M_\varphi$ for $M_{\varphi,\H}$.  It is well-known that if we set
\[
\|\varphi\|_{\mult(\H)}=\|M_\varphi\|_{\mathcal L(\H)},\qquad \varphi\in\mult(\H),
\]
then $\mult(\H)$ becomes a Banach algebra. Moreover, using a standard argument, we have
\begin{equation}\label{eq:vector-propre-multiplier}
M_{\varphi}^*k_\lambda^{\H}=\overline{\varphi(\lambda)}k_\lambda^{\H},\qquad (\lambda\in \Omega),
\end{equation}
which gives
\begin{equation} \label{E:bdd-evaluation-mult}
|\varphi(\lambda)| \leq \|\varphi\|_{\mult(\H)}, \qquad (\lambda \in \Omega).
\end{equation}
See for instance \cite[Chapter 9]{MR3497010} or \cite{Paulsen}.

Let $\H_1,\H_2$ be two RKHS such that $\H_1\subset\H_2$. If $(f_n)_n$ is a sequence in $\H_1$ which is convergent in the weak topology of $\H_2$, we cannot deduce that it also converges in the weak topology of $\H_1$. However, the following result shows that on the bounded subsets of $\H_1$ the above conclusion holds.

\begin{lem}\label{lemA}
Let $\H_1,\H_2$ be two RKHS on a set $\Omega$ such that $\H_1\subset\H_2$, let $(f_n)_n$ be a sequence in  $\H_1$ bounded in $\H_1$-norm by a constant $C$, and let $f \in \H_2$. Assume that $(f_n )_n$ converges to $f$ in the weak topology of $\H_2$. Then the following holds:
\begin{enumerate}[(i)]
\item $f \in \H_1$,
\item $f_n \to f$ in the weak topology of $\H_1$,
\item $\|f\|_{\H_1} \leq C$.
\end{enumerate}
\end{lem}

\proof Since $(f_n)_n$ is uniformly bounded in the norm of $\H_1$, it has a weakly convergent subsequence. More explicitly, there is a subsequence $(f_{n_k})_k$ that converges to some $g\in \H_1$ in the weak topology of $\H_1$. Using \eqref{eq:weakly-convergence-simply-convergence}, we easily see that the two functions $f$ and $g$ coincide on $\Omega$. Therefore $f \in \H_1$. Second, since each $\H_1-$weakly convergent subsequence of $(f_n)_n$ has to converge weakly to $f$ in $\H_1$, we conclude that $(f_n)_n$ itself also converges to $f$ in the weak topology of $\H_1$. Third, the weak convergence in $\H_1$ implies
\[
\|f\|_{\H_1} \leq \liminf_{n \to \infty}\|f_{n}\|_{\H_1} \leq C,
\]
completing the proof. 

\subsection{$H^\infty$ functional calculus for contractions}
Let $T$ be a contraction on a Hilbert space $\H$. We recall that $T$ is said to be {\emph completely non-unitary} if there is no nonzero reducing subspaces $\H_0$ for $T$ such that $T|\H_0$ is a unitary operator. We recall that for a completely non-unitary contraction $T$ on $\H$, we can define an $H^\infty$-functional calculus with the following properties (see \cite[Theorem 2.1, page 117]{Sz-Nagy}):
\begin{enumerate}
\item[(P3)] for every $f\in H^\infty$, we have
\[
\|f(T)\|\leq \|f\|_\infty.
\]
\item[(P4)] If $(f_n)_n$ is a sequence of $H^\infty$ functions which tends boundedly to $f$ on the open unit disc $\mathbb D$ (which means that $\sup_n \|f_n\|_\infty<\infty$ and $f_n(z)\to f(z)$, $n\to+\infty$ for every $z\in\mathbb D$), then $f_n(T)$ tends to $f(T)$ WOT (for the weak operator topology).
\item[(P5)] If $(f_n)_n$ is a sequence of $H^\infty$ functions which tends boundedly to $f$ almost everywhere on $\T=\partial\mathbb D$, then $f_n(T)$ tends to $f(T)$ SOT (for the strong operator topology). 
\end{enumerate}
Finally, we recall that every contraction $T$ on a Hilbert space $\H$ has a unitary dilation $U$ on $\mathcal K$ (which means that $\H\subset\mathcal K$ and $T^n=P_{\H}U^n|\H$, $n\geq 1$) which is minimal (in the sense that $\mathcal K=\bigvee_{-\infty}^{\infty}U^n \H$). 

\subsection{A general framework} \label{S:CATO}

In this note, we consider two analytic reproducing kernel Hilbert spaces $\H_1$ and $\H_2$  on the open unit disc $\mathbb D$ (which means that their elements are analytic on $\mathbb D$) and such that $\H_1\subset \H_2$. A standard application of the closed graph theorem shows that there is a constant $C$ such that
\begin{equation} \label{E:i-bdd}
\|f\|_{\H_2} \leq C \|f\|_{\H_1}, \qquad (f \in \H_1).
\end{equation}
Denote by $\chi$ the function $\chi(z)=z$, $z\in\mathbb D$. Furthermore, we shall assume the following two properties:
\begin{equation}\label{eq:hypothesis1}
\chi\in\mult(\H_2) \quad\mbox{and}\quad \|\chi\|_{\mult(\H_2)}\leq 1,
\end{equation}
and if $X:=M^*_{\chi,\H_2}$ (recall notation \eqref{eq:multiplier-defn}), then
\begin{equation}\label{eq:hypothesis2}
X\H_1\subset \H_1\quad\mbox{and}\quad \|X\|_{\mathcal L(\H_1)}\leq 1.
\end{equation}
The restriction of $X$ to $\H_1$ is denoted by
$$
X_{\H_1}:=X_{|\H_1}.
$$

\subsection{Range Spaces} \label{range}
Let $\mathscr X,\mathscr Y$ be two Hilbert spaces and $T\in\mathcal L(\mathscr X,\mathscr Y)$. We define $\mathcal M(T)$ as the range space equipped with the range norm. More explicitly, $\mathcal M(T)=\mathcal R(T)=T\mathscr X$ and
$$
\|Tx\|_{\mathcal M(T)}=\|P_{(\ker T)^\perp}x\|_{\mathscr X},\qquad x\in\mathscr X,
$$
where $P_{(\ker T)^\perp}$ denotes the orthogonal projection from $\mathscr X$ onto $(\ker T)^\perp$. It is easy to see that $\mathcal M(T)$ is a Hilbert space which is boundedly contained in $\mathscr Y$. A result of Douglas \cite{MR0203464} says that if $A\in\mathcal L(\mathscr X_1,\mathscr Y)$ and $B\in\mathcal L(\mathscr X_2,\mathscr Y)$, then
\begin{equation}\label{eq:douglas1}
\mathcal M(A)\eqcirc \mathcal M(B) \Longleftrightarrow AA^*=BB^*.
\end{equation}
Here the notation $\mathcal M(A)\eqcirc \mathcal M(B)$ means that the Hilbert spaces $\mathcal M(A)$ and $\mathcal M(B)$ coincide as sets and, moreover, have the same Hilbert space structure. We also recall  that if $A,B\in\mathcal L(\mathscr X_1,\mathscr Y)$ and $C\in\mathcal L(\mathscr Y)$, then
\begin{equation}\label{eq:douglas2}
C\mbox{ is a contraction from $\mathcal M(A)$ into $M(B)$} \Longleftrightarrow CAA^*C^*\leq BB^*.
\end{equation}
See also \cite[Corollaries 16.8 and 16.10]{MR3617311}.

\section{Multiplication operators}

Note that \eqref{eq:hypothesis1} implies $\|\chi\|_{\mult(\H_2)}=1$. Indeed, according to \eqref{E:bdd-evaluation-mult}, we have
$$
1=\sup_{z\in\mathbb D}|\chi(z)|\leq \|\chi\|_{\mult(\H_2)}\leq 1.
$$
More generally, since $\bigcap_{n\geq 0}M_{\chi,\mathcal H_2}^n\mathcal H_2=\{0\}$, we see that $M_{\chi,\mathcal H_2}$ is a completely non-unitary contraction. Hence, we get the following consequence.

\begin{lem}\label{lemme-norm-multiplicateur}
Let $\H_2$ be a reproducing kernel Hilbert space of analytic functions on $\mathbb D$ satisfying \eqref{eq:hypothesis1}. Then $\mult(\H_2)=H^\infty$ and for every $\varphi\in H^\infty$, we have $M_{\varphi,\H_2}=\varphi(M_{\chi,\H_2})$ with
\begin{equation}\label{eq:norm-multiplicateur}
\|\varphi\|_{\mult(\H_2)}=\|\varphi\|_\infty.
\end{equation}
\end{lem}

\proof
Let $\varphi\in H^\infty$ and consider the dilates $\varphi_r(z)=\varphi(rz)$, $0<r<1$, $z\in\mathbb D$. If $\varphi(z)=\sum_{n=0}^\infty a_nz^n$ and $f\in \mathcal H_2$, observe that 
$$
\varphi_r(M_{\chi,\H_2})f=\sum_{n=0}^\infty a_n r^n M_{\chi,\H_2}^n f=\sum_{n=0}^\infty a_n r^n \chi^n f=\varphi_r f.
$$
Moreover, by (P5), we have $\varphi_r(M_{\chi,\H_2})f\to \varphi(M_{\chi,\H_2})f$ in $\H_2$ as $r\to 1$. Then, using \eqref{eq:weakly-convergence-simply-convergence}, we get on one hand
$$
\varphi_r(\lambda)f(\lambda)=(\varphi_r(M_{\chi,\H_2})f)(\lambda)\to (\varphi(M_{\chi,\H_2})f)(\lambda),\qquad \mbox{as }r\to 1,\, (\lambda\in\mathbb D),
$$
and on the other hand, $\varphi_r(\lambda)f(\lambda)\to \varphi(\lambda)f(\lambda)$, $r\to 1$ ($\lambda\in\mathbb D$). We thus deduce that $\varphi f=\varphi(M_{\chi,\H_2})f\in\H_2$. In particular, $\varphi\in\mult(\H_2)$ and $M_{\varphi,\H_2}=\varphi(M_{\chi,\H_2})$. Moreover, by (P3), we have 
$$
\|\varphi\|_{\mult(\H_2)}=\|\varphi(M_{\chi,\H_2})\|\leq \|\varphi\|_\infty.
$$
If we combine with \eqref{E:bdd-evaluation-mult}, we get \eqref{eq:norm-multiplicateur}, as claimed.

\begin{lem}\label{Lemme-spectral-measure}
Let $\mathcal H_1$ and $\H_2$ be two reproducing kernel Hilbert spaces of analytic functions on $\mathbb D$ such that $\H_1\subset \H_2$. Assume that $\H_1$ and $\H_2$ satisfy \eqref{eq:hypothesis1} and \eqref{eq:hypothesis2}. Then the minimal unitary dilation of $X_{\H_1}$ has an absolutely continuous scalar spectral measure. In particular, for every $f,g\in\H_1$, there exists $u_{f,g}\in L^1(\T)$ such that
\begin{equation}\label{eq1:spectral-measure}
\langle X_{\H_1}^nf, g \rangle_{\H_1}=\int_{\T}z^n u_{f,g}(z)\,dm(z).
\end{equation}
\end{lem}

\proof
Let $f,g\in\H_1$ and let $\mu_{f,g}$ be the scalar spectral measure associated to the minimal unitary dilation of the contraction $X_{\H_1}$. Then, we have
\begin{equation}\label{eq2:spectral-measure}
\langle X_{\H_1}^nf, g \rangle_{\H_1}=\int_{\T}z^n\,d\mu_{f,g}(z).
\end{equation}
Let us prove that $\mu_{f,g}$ is absolutely continuous with respect to normalized Lebesgue measure $m$ on $\T$. Let $F$ be a closed Borel subset of $\T$ such that $m(F)=0$. Then, we can construct a bounded sequence of polynomials $(q_n)_n$ such that $q_n(z)\to \chi_F(z)$, as $n\to +\infty$, for every $z\in\overline{\mathbb D}$. Indeed, 
Let $f$ be the Fatou function associated to $F$, that is a function $f$ in the disc algebra (that is the closure of polynomials for the sup norm) such that $f=1$ on $F$ and $|f|<1$ on $\overline{\mathbb D}\setminus F$ (See \cite[page 116]{MR2500010} or \cite{Koosis}). Now take $f^n$, $n\geq 0$. The functions $f^n$ are still in the disc algebra. Then if we take $\varepsilon>0$, we can find a polynomial $q_n$ such that 
$$
\sup_{z\in\overline{\mathbb D}}|f^n(z)-q_n(z)|\leq \frac{\varepsilon}{2}.
$$
In particular, we have for every $z\in F$, $|1-q_n(z)|\leq \varepsilon/2$. On the other hand, for $z\in\overline{\mathbb D}\setminus F$, we can find $n_0$ such that for $n\geq n_0$, $|f^n(z)|\leq \varepsilon/2$ (because $|f_n(z)|<1$ and thus $|f^n(z)|\to 0$, as $n\to \infty$). Therefore, for $n\geq n_0$, we have
$$
|q_n(z)|\leq |q_n(z)-f^n(z)|+|f^n(z)|\leq \frac{\varepsilon}{2}+\frac{\varepsilon}{2}=\varepsilon.
$$ 
Hence $q_n(z)$ tends to $1$ for $z\in F$ and to $0$ for $z\in\overline{\mathbb D}\setminus F$. In other words, $q_n$ tends  to $\chi_F$ pointwise. On the other hand, we have of course
$$
\sup_{z\in\mathbb D}|q_n(z)|\leq 1+\frac{\varepsilon}{2},
$$
which proves that the sequence $(q_n)_n$ is also bounded, and we are done.

Now, since $(q_n)_n$ converges boundedly to $0$ on $\mathbb D$ and since $X$ is a completely unitary contraction, we deduce from (P4) that $(q_n(X))_n$ converges WOT to $0$ in $\mathcal L(\H_2)$. Hence it implies that $(q_n(X_{\H_1})f)_n$ converges weakly to $0$ in $\H_2$.  On the other hand, by von Neumann inequality, we have
$$
\|q_n(X_{\H_1})f\|_{\H_1}\leq \|q_n\|_{\infty}\|f\|_{\H_1}\leq C \|f\|_{\H_1},
$$
where $C=\sup_{n}\|q_n\|_{\infty}<+\infty$. By Lemma~\ref{lemA}, we deduce that $(q_n(X_{\H_1})f)_n$ converges weakly to $0$ in $\H_1$. But, according to \eqref{eq2:spectral-measure}, we have 
$$
\langle q_n(X_{\H_1})f, g \rangle_{\H_1}=\int_{\T}q_n(z)\,d\mu_{f,g}(z),
$$
which gives that
$$
\lim_{n\to+\infty}\int_{\T}q_n(z)\,d\mu_{f,g}(z)=0.
$$
It remains to apply dominated Lebesgue convergence theorem to get
$$
\int_{\T}\chi_F(z)\,d\mu_{f,g}(z)=0,
$$
which implies that $\mu_{f,g}(F)=0$. Hence $\mu_{f,g}$ is absolutely continuous with respect to $m$, as claimed.
\begin{thm}\label{C} Let $\H_1$ and $\H_2$  be two reproducing kernel Hilbert spaces of analytic functions on $\mathbb D$ such that $\H_1\subset \H_2$. Assume that $\H_1$ and $\H_2$ satisfy \eqref{eq:hypothesis1} and \eqref{eq:hypothesis2}. Let $\varphi \in H^\infty$.  Then  $M_{\varphi,\H_2}^*$ maps $\H_1$ into itself, and if $f,g\in\H_1$, we have
\begin{equation}\label{eq:spectral-Mvarphi}
\langle M_{\varphi,\H_2}^*f, g \rangle_{\H_1}=\int_{\T}\varphi^*(z)u_{f,g}(z)\,dm(z),
\end{equation}
where $\varphi^*(z)=\overline{\varphi(\overline{z})}$. 
\end{thm}

\proof
Let us first assume that $\varphi$ is holomorphic on $\overline{\mathbb D}$ and let us consider the Taylor series of $\varphi$,  $\varphi(z)=\sum_{n=0}^\infty a_n z^n$. Then we have 
\begin{equation}\label{eq:calcul-fonctionnel-holomorphe-basique}
M_{\varphi,\H_2}^*=\varphi(M_{\chi,\H_2})^*=\sum_{n=0}^\infty \overline{a_n}X^n.
\end{equation}
Since $X\H_1\subset\H_1$, the last equation implies that $M_{\varphi,\H_2}^*\H_1\subset\H_1$. Now using that $\sum_{n=0}^{\infty}|a_n|<\infty$ and \eqref{eq1:spectral-measure}, we get 
\begin{eqnarray*}
\langle M_{\varphi,\H_2}^*f, g \rangle_{\H_1}&=& \sum_{n=0}^\infty \overline{a_n}\langle X_{\H_1}^n f,g\rangle_{\H_1} \\
&=&\sum_{n=0}^\infty \overline{a_n} \int_{\mathbb T}z^n u_{f,g}(z)\,dm(z)\\
&=&\int_{\mathbb T}\sum_{n=0}^\infty \overline{a_n} z^n u_{f,g}(z)\,dm(z)\\
&=&\int_{\mathbb T}\varphi^*(z) u_{f,g}(z)\,dm(z).
\end{eqnarray*}
This proves \eqref{eq:spectral-Mvarphi} for $\varphi$ which is holomorphic on $\overline{\mathbb D}$. We also observe that
\begin{eqnarray*}
\left| \langle M_{\varphi,\H_2}^*f, g \rangle_{\H_1}\right| &\leq& \int_{\mathbb T}|\varphi^*(z)| |u_{f,g}(z)|\,dm(z)\\
&\leq & \|\varphi\|_{\infty} \int_{\mathbb T}|u_{f,g}(z)|\,dm(z).
\end{eqnarray*}
But by spectral theorem, we know that $\int_{\mathbb T}|u_{f,g}(z)|\,dm(z)=\|\mu_{f,g}\|\leq \|f\|_{\H_1}\|g\|_{\H_1}$, which gives 
\begin{equation}\label{eq:generalisation-v-n-algebre-disque}
\|M_{\varphi,\H_2}^*f\|_{\H_1}\leq \|\varphi\|_{\infty}\|f\|_{\H_1}.
\end{equation}
 Now let $\varphi\in H^\infty$ and define the dilates $\varphi_r(z)=\varphi(rz)$, $0<r<1$, $z\in\mathbb D$. Observe that $\varphi_r$ are holomorphic on $\overline{\mathbb D}$. By the previous argument, we get that $M_{\varphi_r,\H_2}^*$ maps $\H_1$ into itself and 
\begin{equation}\label{eq:formule-pour-varphir}
\langle M_{\varphi_r,\H_2}^*f, g \rangle_{\H_1}=\int_{\mathbb T}\varphi_r^* u_{f,g}\,dm,\qquad f,g\in\H_1.
\end{equation}
Since $\varphi_r$ converges boundedly to $\varphi$ on $\mathbb D$ as $r\to 1$, and since $M_{\chi,\H_2}$ is a completely non unitary contraction on $\H_2$, we get that $M_{\varphi_r,\H_2}^*f$ converges weakly to $M_{\varphi,\H_2}^*f$ in $\H_2$ as $r\to 1$. On the other hand, using \eqref{eq:generalisation-v-n-algebre-disque}, we have
$$
\|M_{\varphi_r,\H_2}^*f\|_{\H_1}\leq \|\varphi_r\|_{\infty} \|f\|_{\H_1}\leq \|\varphi\|_{\infty} \|f\|_{\H_1}.
$$
Lemma~\ref{lemA} now implies that $M_{\varphi,\H_2}^*f$ belongs to $\H_1$ and $M_{\varphi_r,\H_2}^*f$ converges weakly to $M_{\varphi,\H_2}^*f$ in $\H_1$ as $r\to 1$. Letting $r\to 1$ in \eqref{eq:formule-pour-varphir} and using dominated convergence, we deduce that formula \eqref{eq:spectral-Mvarphi} is satisfied by $\varphi$, completing the proof.

\begin{rem}\label{remark-multiplicateurnorme-on-H1}
It follows immediately from \eqref{eq:spectral-Mvarphi} that for $\varphi\in H^\infty$, we have
$$
\|M_{\varphi,\H_2}^*\|_{\mathcal L(\H_1)}\leq \|\varphi\|_{\infty}.
$$
\end{rem}

Given a bounded operator $T$ on a Hilbert space $\H$, the family of all closed $T$-invariant subspaces of $\H$ is denoted by $\mbox{Lat}(T)$.

\begin{cor}\label{cor:invariant}
Let $\H_1$ and $\H_2$  be two reproducing kernel Hilbert spaces of analytic functions on $\mathbb D$ such that $\H_1\subset \H_2$. Assume that $\H_1$ and $\H_2$ satisfy \eqref{eq:hypothesis1} and \eqref{eq:hypothesis2}. Then, for every $\varphi\in H^\infty$, we have
$$
\mbox{Lat}(X_{\H_1})\subset \mbox{Lat}(M_{\varphi,\H_2}^*|\H_1).
$$
\end{cor}

\proof
Let $\varphi\in H^\infty$, $\varphi_r(z)=\varphi(rz)$, $0<r<1$, and let $E\in \mbox{Lat}(X_{\H_1})$. Note that \eqref{eq:calcul-fonctionnel-holomorphe-basique} implies that $M_{\varphi_r,\H_2}^* E\subset E$. On the other hand, as we have seen in the proof of Theorem~\ref{C}, $M_{\varphi_r,\H_2}^* \to M_{\varphi,\H_2}^*$, as $r\to 1$, in the weak operator topology of $\mathcal{L}(\H_1)$. Since a norm-closed subspace is also weakly closed \cite{MR1892228}, we conclude that $M_{\varphi,\H_2}^* E \subset E$, as claimed.

To conclude this section, we show that Theorem~\ref{C} has an interesting application in relation with the F-property. Recall that a linear manifold $V$ of $H^1$ is said to have the F-property if whenever $f \in V$ and $\theta$ is an inner function which is lurking in $f$, i.e., $f/\theta \in H^1$ or equivalently $\theta$ divides the inner part of $f$, then we actually have $f/\theta \in V$. This concept was first introduced by V. P. Havin \cite{MR0289783} and it plays a vital role in the analytic function space theory. Several classical spaces have the F-properties. the list includes Hardy spaces $H^p$, Dirichlet space $\mathcal{D}$, BMOA, VMOA, and the disc algebra $\mathcal{A}$. See \cite{MR2032687, MR2290751,Shirokov}. However, for the Bloch spaces $\mathfrak{B}$ and $\mathfrak{B}_0$, we know that $\mathfrak{B} \cap H^p$ and $\mathfrak{B}_0 \cap H^p$ do not have the F-property \cite{MR2199173}. Using the tools developed in Section \ref{S:CATO}, we will see that in the situation when $\H_1\subset H^2$ satisfies \eqref{eq:hypothesis2}, then $\H_1$ has the F-property. First, let us note that  $\H_2=H^2$ satisfies \eqref{eq:hypothesis1} and $M_{\chi,\H_2}=S$ is the classical forward shift operator. Thus, $X=M_{\chi,\H_2}^*=S^*$ is the backward shift operator
$$
(S^*f)(z)=\frac{f(z)-f(0)}{z},\qquad f\in H^2,z\in\mathbb D.
$$
In this context, if $\H_1$ is a reproducing kernel Hilbert space such that $\H_1\subset H^2$, the condition \eqref{eq:hypothesis2} can be rephrased as
\begin{equation}\label{eq:hypothese-H1}
S^*\H_1\subset \H_1\quad\mbox{and }\quad \|S^*|\H_1\|\leq 1.
\end{equation}
Recall that for $\psi\in L^\infty(\mathbb T)$, the Toeplitz operator $T_\psi$ is defined  on $H^2$ by $T_\psi(f)=P_{+}(\psi f)$ where $P_+$ is the Riesz projection (the orthogonal projection from $L^2(\mathbb T)$ onto $H^2$). If $\varphi\in H^\infty=\mult(H^2)$, then $M_{\varphi,H^2}=T_\varphi$ and $M_{\varphi,H^2}^*=T_{\overline\varphi}$. In this situation, we get the following result.

\begin{thm}\label{T:F-property}
Let $\H_1$ be a reproducing kernel Hilbert space contained in $H^2$, and assume that it satisfies \eqref{eq:hypothese-H1}. Then the space $\H_1$ has the F-property. Moreover, if $f\in \H_1$ and $\theta$ is an inner function which divides $f$, then
\[
\left\| \frac{f}{\theta} \right\|_{\H_1} \leq \|f\|_{\H_1}.
\]
\end{thm}

\proof
Assume that $f\in\H_1$ and that $\theta$ is an inner function so that $f/\theta\in H^1$. In fact, by Smirnov Theorem \cite{MR2500010}, we actually have $\psi:=f/\theta \in H^2$. Therefore, 
\[
T_{\overline{\theta}} (f) = P_+(\overline{\theta} f)= P_+(\psi) = \psi.
\]
But according to Theorem~\ref{C} and Remark~\ref{remark-multiplicateurnorme-on-H1}, $T_{\overline\theta}$ acts contractively on $\H_1$. Hence $\psi=T_{\overline{\theta}} (f)\in\H_1$ and 
\[
\left\| \frac{f}{\theta} \right\|_{\H_1}=\left\|T_{\overline\theta}f\right\|_{\H_1} \leq \|f\|_{\H_1},
\]
as claimed. 

\section{Invariant subspaces and cyclicity}
The following result says that under certain circumstances, the closed invariant subspaces of $X_{\H_1}=X_{|\H_1}$ are exactly the trace on $\H_1$ of the closed invariant subspaces of $X$. Despite the following characterization, the implication $(i) \Longrightarrow (ii)$ is the essential part of the result.

\begin{thm}\label{Thm:invariant-subspace-backward}
Let $\H_1$ and $\H_2$  be two analytic reproducing kernel Hilbert spaces on $\mathbb D$ such that $\H_1\subset \H_2$ and satisfying \eqref{eq:hypothesis1} and \eqref{eq:hypothesis2}. Assume that there exists an outer function $\varphi\in H^\infty$ such that $\mathcal R(M_{\varphi,\H_2}^*)\subset\H_1$. Then, for every $\mathcal E\subset\H_1$, the following assertions are equivalent.
\begin{enumerate}[(i)]
\item $\mathcal E$ is a closed subspace of $\H_1$ invariant under $X_{\H_1}$;
\item there is a closed subspace $E$ of $\H_2$ invariant under $X=M_{\chi,\H_2}^*$ such that $\mathcal E=E\cap \H_1$.
\end{enumerate}
Moreover, $\mathcal E=\H_1$ if and only if $E=\H_2$.
\end{thm}
The proof will be based on the following lemma, which extends \cite[Lemmata 17.21 and 24.30]{MR3617311} in our general context.

\begin{lem}\label{lem-densite}
Let $\H_1$ and $\H_2$  be two analytic reproducing kernel Hilbert spaces on $\mathbb D$ such that $\H_1\subset \H_2$ and satisfying \eqref{eq:hypothesis1} and \eqref{eq:hypothesis2}. Assume that there exists an outer function $\varphi\in H^\infty$ such that $\mathcal R(M_{\varphi,\H_2}^*)\subset\H_1$. Then, for every $\mathcal E\in\mbox{Lat}(X_{\H_1})$, the space $M_{\varphi,\H_2}^*\mathcal E$ is dense in $\mathcal E$ with respect to the norm topology of $\H_1$.
\end{lem}

\proof
According to Corollary~\ref{cor:invariant}, we know that $M_{\varphi,\H_2}^*\mathcal E\subset\mathcal E$. Now let $g\in\mathcal E$, $g\perp M_{\varphi,\H_2}^*\mathcal E$ in the $\H_1$-topology. In particular, for every $n\geq 0$, we have
$$
0=\langle M_{\varphi,\H_2}^*X_{\H_1}^n g,g \rangle_{\H_1}.
$$
Observe now that $M_{\chi,\H_2}^n M_{\varphi,\H_2}=M_{\chi^n\varphi,\H_2}$, which gives $M_{\varphi,\H_2}^*X^n=M_{\chi^n\varphi,\H_2}^*$. Hence, by Theorem~\ref{C}, we get
$$0=\langle M_{\chi^n\varphi,\H_2}^* g,g \rangle_{\H_1}=\int_{\mathbb T}\varphi^*(z) z^n u_{g,g}(z)\,dm(z),
$$
for every $n\geq 0$. We thus deduce that $\varphi^* u_{g,g}\in H_0^1$. Since $\varphi^*$ is outer and $u_{g,g}\in L^1(\mathbb T)$, Smirnov Theorem \cite{Duren} implies that $u_{g,g}\in H_0^1$. Since $u_{g,g}\geq 0$, this gives $u_{g,g}=0$, that is $g=0$, completing the proof of the Lemma. 

\noindent
\emph{Proof of Theorem~\ref{Thm:invariant-subspace-backward}.} \hskip 0,1cm $(ii)\Longrightarrow (i)$: Let $E$ be a closed subspace of $\H_2$, invariant under $X=M_{\chi,\H_2}^*$ such that $\mathcal E=E\cap \H_1$. First, let us check that $\mathcal E$ is a closed subspace of $\H_1$.  The verification essentially owes to \eqref{E:i-bdd}. To do so, let $f \in \H_1$ be in the $\H_1$-closure of $E \cap \H_1$. Then there is a sequence $(f_n)_n$ in $E\cap \H_1$ which converges to $f$ in the norm topology of $\H_1$. Since $\H_1$ is boundedly contained in $\H_2$, the sequence $(f_n)_n$ also converges to $f$  in $\H_2$. Since $E$ is closed in $\H_2$, the function $f$ must belong to $E$. Hence, $f\in \mathcal E=E\cap \H_1$, which proves that  $\mathcal E$ is closed in $\H_1$. The fact that $\mathcal E$ is invariant under $X_{\H_1}=M_{\chi,\H_2}^*|\H_1$ is immediate.  \\

$(i)\Longrightarrow (ii)$: A standard argument using the closed graph theorem implies that, according to $\mathcal R(M_{\varphi,\H_2}^*)\subset\H_1$, the mapping $M_{\varphi,\H_2}^*$ from $\H_2$ into $\H_1$ is a bounded operator. Now let $\mathcal E$ be a closed subspace of $\H_1$, and assume that $\mathcal E$ is invariant under $X_{\H_1}$. Denote by $E$ the closure of $\mathcal E$ in the $\H_2$-topology. It is clear that $E$ is a closed subspace of $\H_2$ which is invariant under $X$. Let us prove that $\mathcal E=E\cap\H_1$.  The inclusion $\mathcal E\subset E\cap\H_1$ is trivial. For the reverse inclusion, let us verify that
\begin{equation}\label{eq:Beurling-type-thm1}
M_{\varphi,\H_2}^* E\subset\mathcal E.
\end{equation}
Let $f\in E$. By definition, there is a sequence $(f_n)_n$  in $\mathcal E$ which converges to $f$ in the $\H_2$-topology. Then, since $M_{\varphi,\H_2}^*$ is bounded from $\H_2$ into
$\H_1$, the sequence $(M_{\varphi,\H_2}^*f_n)_n$ tends to $M_{\varphi,\H_2}^*f$ in the $\H_1$-topology. Since $f_n\in\mathcal E$, Corollary~\ref{cor:invariant} implies that $M_{\varphi,\H_2}^*f_n\in\mathcal E$ and since $\mathcal E$ is closed in $\H_1$, then $M_{\varphi,\H_2}^*f\in\mathcal E$, which proves \eqref{eq:Beurling-type-thm1}. In particular, we have
$$
M_{\varphi,\H_2}^*(E\cap\H_1)\subset \mathcal E,
$$
and since $E\cap \H_1$ is a closed subspace of $\H_1$ invariant with respect to $X_{\H_1}$, it follows from Lemma~\ref{lem-densite} that $M_{\varphi,\H_2}^*(E\cap\H_1)$ is dense in $E\cap\H_1$, which implies
$$
E\cap\H_1\subset \mathcal E.
$$
Thus we have $\mathcal E=E\cap \H_1$. \\

It remains to prove that $\mathcal E=\H_1$ if and only if $E=\H_2$. Assume first that $\mathcal E=\H_1$. Then $\mathcal R(M_{\varphi,\H_2}^*)\subset E$. But note that $\ker\,(M_{\varphi,\H_2})=\{0\}$, whence $\mathcal R(M_{\varphi,\H_2}^*)$ is dense in $\H_2$. Hence we get $E=\H_2$. Conversely, assume that $E=\H_2$. Then
$$
\mathcal E=E\cap\H_1=\H_2\cap\H_1=\H_1,
$$
which concludes the proof. 

\begin{rem}
As already noted, $\mathcal R(M_{\varphi,\H_2}^*)$ is always  dense in $\H_2$ and thus under the hypothesis of Theorem~\ref{Thm:invariant-subspace-backward} (that is if there exists an (outer) function $\varphi\in H^\infty$ such that $\mathcal R(M_{\varphi,\H_2}^*)\subset\H_1$), then automatically $\H_1$ is dense in $\H_2$. 
\end{rem}

Theorem~\ref{Thm:invariant-subspace-backward} has an immediate application in characterization of cyclic vectors.
\begin{cor}\label{Cor:cyclic-general}
Let $\H_1$ and $\H_2$  be two analytic reproducing kernel Hilbert spaces on $\mathbb D$ such that $\H_1\subset \H_2$ and satisfying \eqref{eq:hypothesis1} and \eqref{eq:hypothesis2}. Suppose that there exists an outer function $\varphi\in H^\infty$ such that $\mathcal R(M_{\varphi,\H_2}^*)\subset\H_1$.
Let $f\in\H_1$. Then the following assertions are equivalent.
\begin{enumerate}[(i)]
\item $f$ is cyclic for $X=M_{\chi,\H_2}^*$.
\item $f$ is cyclic for $X_{\H_1}=X_{|\H_1}$.
\end{enumerate}
\end{cor}

\proof
$(i)\Longrightarrow (ii)$: Assume that $f$ is cyclic for $X$ in $\H_2$ and denote by $\mathcal E$ the subspace of $\H_1$ defined by
$$
\mathcal E=\overline{\mbox{Span}(X_{\H_1}^nf:n\geq 0)}^{\H_1}.
$$
It is clear that $\mathcal E$ is a closed subspace of $\H_1$, invariant with respect to $X_{\H_1}$. Assume that $\mathcal E\neq \H_1$. Then, according to Theorem~\ref{Thm:invariant-subspace-backward},  there exists a closed subspace $E$ of $\H_2$, $E\neq \H_2$, invariant under $X$ such that $\mathcal E=E\cap \H_1$. In particular, $f\in E$, and thus it is not cyclic for $X$, which is contrary to the hypothesis. Thus $\mathcal E=\H_1$ and $f$ is cyclic for $X_{\H_1}$.

$(ii)\Longrightarrow (i)$: Let $g\in\mathcal R(M_{\varphi,\H_2}^*)$. Then  $g\in\H_1$ and if $f$ is cyclic for $X_{\H_1}$, there exists a sequence of polynomials $(p_n)$ such that
$$
\|p_n(X_{\H_1})f-g\|_{\H_1}\to 0,\qquad \mbox{as }n\to\infty.
$$
Since $\H_1$ is contained boundedly in $\H_2$, then we have
$$
\|p_n(X_{\H_2})f-g\|_{\H_2}\to 0,\qquad \mbox{as }n\to\infty.
$$
Thus,  $\mathcal R(M_{\varphi,\H_2}^*)\subset \overline{\mbox{Span}({X^*}^nf:n\geq 0)}^{\H_2}$. Since $\mathcal R(M_{\varphi,\H_2}^*)$ is dense in $\H_2$, we get 
$$
\overline{\mbox{Span}({X^*}^nf:n\geq 0)}^{\H_2}=\H_2,
$$
completing the proof. 

We can apply Theorem~\ref{Thm:invariant-subspace-backward}  and Corollary~\ref{Cor:cyclic-general} to some specific reproducing kernel Hilbert spaces contained in the Hardy space $H^2$ on $\mathbb D$.

\begin{thm}\label{Thm:invariant-subspace-backwardbis}
Let $\H_1$  be a reproducing kernel Hilbert space contained in $H^2$ that satisfies \eqref{eq:hypothese-H1} and assume that there exists an outer function $\varphi\in H^\infty$ such that $T_{\overline\varphi}H^2\subset\H_1$.
Then, for every $\mathcal E\subsetneq\H_1$, the following assertions are equivalent.
\begin{enumerate}[(i)]
\item $\mathcal E$ is a closed subspace of $\H_1$ invariant under $X_{\H_1}$;
\item there is an inner function $\Theta$ such that $\mathcal E=K_\Theta \cap \H_1$.
\end{enumerate}
Moreover, if $f\in\H_1$, then $f$ is cyclic for $S^*|\H_1$ if and only if $f$ is cyclic for $S^*$.
\end{thm}
\proof
It is sufficient to combine Theorem~\ref{Thm:invariant-subspace-backward} and  Corollary~\ref{Cor:cyclic-general} with Beurling's theorem.

\begin{rem}
Note that the hypothesis $T_{\overline\varphi}H^2\subset\H_1$ implies that polynomials belong to $\H_1$.
\end{rem}

Regarding the last part of Theorem~\ref{Thm:invariant-subspace-backwardbis}, let us mention that a well-known theorem of Douglas--Shapiro--Shields \cite{MR0270196} says that a function $f$ in $H^2$ is cyclic for $S^*$ if and only if $f$ has no bounded type meromorphic pseudo continuation across $\mathbb T$ to $\mathbb{D}_e=\{z:1<|z|\leq \infty\}$.

\section{Applications}
In this section, we give some examples of RKHS for which our main Theorem \ref{Thm:invariant-subspace-backward} can be applied.

\subsection {A general RKHS}\label{subsection-generalRKHS}
Let $\H_2$ be an analytic RKHS on $\mathbb D$ satisfying \eqref{eq:hypothesis1}. Let $\varphi \in H^\infty$ and $\H_1:= \MM(M^*_{\varphi, \H_2})$. Recall the definition of the range space from Section~\ref{range}. Then $\H_1$ is also an analytic RKHS on $\mathbb D$ which is contained in $\H_2$. Observe that $\H_1$ satisfies \eqref{eq:hypothesis2}. Indeed, since $M_{\varphi, \H_2}M_{\chi, \H_2}=M_{\chi, \H_2}M_{\varphi, \H_2}$, we have $XM^*_{\varphi, \H_2}=M^*_{\varphi, \H_2}X$, which implies that $X\H_1 \subset \H_1$. Moreover, if $f = M^*_{\varphi,\H_2}g \in \H_1$ for some $g\in (\ker M^*_{\varphi,\H_2})^{\perp}$, then
\begin{eqnarray*}
\|Xf\|_{\H_1}&=& \|X M^*_{\varphi,\H_2}g\|_{\H_1}\\
&=& \| M^*_{\varphi,\H_2} X g\|_{\H_1}\\
&=& \|Xg\|_{\H_2}\\
&\leq& \|g\|_{\H_2}=\|f\|_{\H_1}
\end{eqnarray*}
Thus $\HH_1$ satisfies (\ref{eq:hypothesis2}). In this context, we get immediately from Theorem~\ref{Thm:invariant-subspace-backward} the following.
\begin{cor}\label{cor-generalRKHS-range-space}
Let $\H_2$ be an analytic RKHS satisfying \eqref{eq:hypothesis1}. Let $\varphi$ be an outer function in $H^\infty$ and let $\H_1:= \MM(M^*_{\varphi, \H_2})$. Then, for every $\mathcal E\subset\H_1$, the following assertions are equivalent.
\begin{enumerate}[(i)]
\item $\mathcal E$ is a closed subspace of $\H_1$, invariant under $X_{\H_1}$;
\item There is a closed subspace $E$ of $\H_2$, invariant under $X=M_{\chi,\H_2}^*$ such that $\mathcal E=E\cap \H_1$.
\end{enumerate}
Moreover $\mathcal E=\H_1$ if and only if $E=\H_2$.
\end{cor}

\subsection{The space $\mathcal M(\overline\varphi)$}
Let $\H_2=H^2$ be the Hardy space on $\mathbb D$, $\varphi$ an outer function in $H^\infty$ and $\H_1=\mathcal M(T_{\overline\varphi})$ which we denote for simplicity $\mathcal M(\overline\varphi)$. The space $H^2$ trivially satisfies \eqref{eq:hypothesis1}  and according to the discussion at the beginning of
Subsection~\ref{subsection-generalRKHS}, the space $\mathcal M(\overline{\varphi})$ is an analytic RKHS contained in $H^2$ and satisfying  \eqref{eq:hypothesis2} (or equivalently \eqref{eq:hypothese-H1}). Again, for simplicity, we write  $X_{\overline\varphi}=X_{\mathcal M(\overline\varphi)}=S^*|\mathcal M(\overline\varphi)$.

In this context, we can apply Theorem~\ref{Thm:invariant-subspace-backwardbis} which immediately gives the following. 
\begin{cor}\label{Cor:Beurling}
Let $\varphi$ be an outer function. Then the following assertions are equivalent.
\begin{enumerate}[(i)]
\item $\mathcal E$ is a closed subspace of $\mathcal M(\overline\varphi)$,
$\mathcal E\neq\mathcal M(\overline\varphi)$, and $\mathcal E$ is invariant under $X_{\overline\varphi}$.
\item There is an inner function $\Theta$ such that $\mathcal E=K_\Theta\cap \mathcal M(\overline\varphi)$.
\end{enumerate}
\end{cor}

\subsection{de Branges--Rovnyak space $\H(b)$}
Let $b\in \mathbf b(H^\infty)$--the closed unit ball of $H^\infty$. The de Branges--Rovnyak space $\H(b)$ is defined as
$$
\H(b)=\mathcal M((I-T_bT_{\overline b})^{1/2}).
$$
For details on de Branges--Rovnyak spaces, we refer to \cite{MR3617311,sarason1994sub}. Here we will just recall what will be useful for us. 

It is well--known that $\H(b)$ is an analytic RKHS contractively contained in $H^2$ and invariant with respect to $S^*$. Moreover, the operator $X_b=S^{*}|\H(b)$ acts as a contraction on $\H(b)$. In particular, the space $\H(b)$ satisfies the hypothesis \eqref{eq:hypothese-H1}. Assume now that $b$ is a non-extreme point of $\mathbf b(H^\infty)$, meaning that $\log(1-|b|)\in L^1(\T)$. Thus, there exists a unique outer function $a$ such that $a(0)>0$ and $|a|^2+|b|^2=1$ a.e. on $\mathbb T$. This function $a$ is called the \emph{pythagorean mate} of $b$. It is well-known that $\mathcal R(T_{\overline a})\subset \H(b)$.

We can then apply Theorem~\ref{Thm:invariant-subspace-backwardbis} to $\H_1=\H(b)$ and $\H_2=H^2$ to recover the following result due to Sarason (\cite{Sarason-86-OT}, Theorem 5).
\begin{cor}[Sarason]\label{thm:Sarason-backward}
Let $b$ be a non-extreme point of the closed unit ball of $H^\infty$, and let $\mathcal E$ be a closed subspace of $\H(b)$, $\mathcal E\neq\H(b)$. Then the following are equivalent.
\begin{enumerate}[(i)]
\item $\mathcal E$ is invariant under $X_b$.
\item There exists an inner function $\Theta$ such that $\mathcal E=K_\Theta\cap\H(b)$.
\end{enumerate}

\end{cor}

\begin{rem}
As already noted, hypothesis of Theorem~\ref{Thm:invariant-subspace-backwardbis} implies that polynomials belongs to $\H_1$. In the case when $\H_1=\H(b)$, we know that it necessarily implies that $b$ is non-extreme. In the extreme case, the backward shift invariant subspaces have been described by Suarez \cite{Suarez-Indiana-97}, also using some Sz.-Nagy-Foias model theory, but  the situation is rather more complicated.
\end{rem}

\subsection{Sub-Bergman Hilbert space}
The Bergman space $A^2$ on $\mathbb D$ is defined as the space of analytic functions $f$ on $\mathbb D$ satisfying
$$
\|f\|_{A^2}^2:=\int_\mathbb D |f(z)|^2 dA(z) < \infty,
$$
where $dA(z)$ is the normalized area measure on $\mathbb D$. In \cite{sultanic2006sub,zhu1996sub,MR1990528}, an analogue of de Branges--Rovnyak spaces was considered in this context. Recall that the Toeplitz operator on $A^2(\mathbb D)$ with symbol $\varphi\in L^\infty(\mathbb D)$ is defined as
$$ T_{\varphi}(f)=P_{A^2}(\varphi f),$$
where $P_{A^2}$ is the Bergman projection (that is the orthogonal projection from $L^2(\mathbb D,dA)$ onto $A^2$).  It is clear that $T_\varphi^*=T_{\overline\varphi}$. Given $\varphi\in L^\infty(\mathbb D)$, we define the sub--Bergman Hilbert space $\H(\varphi)$ as
$$
\H(\varphi)=\mathcal M((I-T_\varphi T_{\overline\varphi})^{1/2}).
$$
In other words, $\H(\varphi)=(I-T_\varphi T_\varphi^*)^{1/2}A^2$ and it is equipped with the inner product
$$\langle (I-T_\varphi T_\varphi^*)^{1/2}f,(I-T_\varphi T_\varphi^*)^{1/2}g \rangle_{\HH(\varphi)}:=\langle f, g\rangle_{A^2},$$
for every $f,g\in A^2\ominus\ker (I-T_\varphi T_\varphi^*)$. We keep the same notation as the de Branges--Rovnyak spaces, but there will be no ambiguity because in this subsection,  the ambient space is $A^2$ (in contrast with the de Branges--Rovnyak spaces for which the ambient space is $H^2$). We refer the reader to \cite{zhu1996sub} for details about this space.

The shift operator (also denoted $S$ in this context), defined as $S=T_z$, is clearly a contraction and $S^*=T_{\overline z}$. As we have seen, the de Branges--Rovnyak spaces are invariant with respect to the backward shift operator which acts as a contraction on them. In the context of sub--Bergman Hilbert spaces, the analogue of this property is also true. The proof is the same but we include it for completeness.
\begin{lem}\label{Lem:Hb-Sstar-invariant}
Let $b\in \mathbf b(H^\infty)$. Then $S^*$ acts as a contraction on $\H(b)$.
\end{lem}
\proof
We first prove that $S^*$ acts as a contraction on $\H(\overline b)$. According to \eqref{eq:douglas2}, we should prove that
$$
S^*(I-T_{\overline b}T_b)S\leq I-T_{\overline b}T_b,
$$
that is
\begin{equation}\label{eq:Lem-Hb-Sstar-invariant}
T_{\overline z}(I-T_{\overline b}T_b)T_z\leq I-T_{\overline b}T_b.
\end{equation}
But, if $\varphi,\psi\in L^\infty(\mathbb D,dA)$ and at least one of them is in $H^\infty$, then
\begin{equation}\label{eq2:Lem-Hb-Sstar-invariant}
T_{\overline\psi}T_\varphi=T_{\overline{\psi}\varphi}.
\end{equation}
See \cite[Proposition 7.1]{zhu2007operator}. Then \eqref{eq:Lem-Hb-Sstar-invariant} is equivalent to $T_{|z|^2(1-|b|^2)}\leq T_{1-|b|^2}$, that is
$$
0\leq T_{(1-|z|^2)(1-|b|^2)}.
$$
Since $(1-|z|^2)(1-|b|^2) \geq 0$ on $\mathbb D$, the latter inequality is satisfied (see also \cite[Proposition 7.1]{zhu2007operator}) and thus $S^*$ acts as a contraction on $\H(\overline b)$.

To pass to the $\H(b)$ case, we use a well--known relation between $\H(b)$ and $\H(\overline b)$: let $f\in A^2$; then $f\in\H(b)$ if and only if $T_{\overline b}f\in\H(\overline b)$ and
$$
\|f\|_{\H(b)}^2=\|f\|_{A^2}^2+\|T_{\overline b}f\|^2_{\H(\overline b)}.
$$
So let $f\in\H(b)$. Since $T_{\overline b}S^* f=S^* T_{\overline b}f$ and $\H(\overline b)$ is invariant with respect to $S^*$, we get that $T_{\overline b}S^*f\in \H(\overline b)$, whence $S^*f\in\H(b)$ and
\begin{eqnarray*}
\|S^*f\|_{\H(b)}^2&=&\|S^* f\|_{A^2}^2+\|T_{\overline b}S^* f\|_{\H(\overline b)}^2\\
&=& \|S^* f\|_{A^2}^2+\|S^*T_{\overline b} f\|_{\H(\overline b)}^2\\
&\leq & \| f\|_{A^2}^2+\|T_{\overline b} f\|_{\H(\overline b)}^2=\|f\|_{\H(b)}^2.
\end{eqnarray*}
Hence $S^*$ is a contraction on $\H(b)$, completing the proof.

\smallskip
According to Lemma~\ref{Lem:Hb-Sstar-invariant}, we see that $A^2$ satisfies \eqref{eq:hypothesis1} and $\H(b)$ satisfies \eqref{eq:hypothesis2}. We will show that under the additional hypothesis that $b$ is a non-extreme point of the closed unit ball of $H^\infty$, we can apply our Corollary~\ref{cor-generalRKHS-range-space} to $\H_1=\H(b)$ and $\H_2=A^2$.

\begin{cor}\label{Cor:sub-bergman-space}
Let $b$ be a non--extreme point of the unit ball of $H^\infty$ and $a$ its pythagorean mate. Then the following are equivalent.
\begin{enumerate}[(i)]
\item $\mathcal E$ is a closed subspace of $\H(b)$, invariant under $X_b=S^*|\H(b)$.
\item There is a closed subspace $E$ of $A^2$, invariant under $S^*$, such that $\mathcal E=E\cap \H(b)$.
\end{enumerate}
\end{cor}
\proof
It is known that since $b$ is analytic, then $\H(b)=\H(\bar b)$ with equivalent norms, see \cite[page 641]{sultanic2006sub}. Moreover, according to \eqref{eq2:Lem-Hb-Sstar-invariant}, we have
$$
I-T_{\bar b}T_b=T_{1-|b|^2}=T_{|a|^2}=T_{\bar a}T_a.
$$
This identity implies by \eqref{eq:douglas1} that $\H(\bar b)\eqcirc\MM(T_{\bar a})$. Hence $\H(b)=\MM(T_{\bar a})=\MM(T_a^*)$ with equivalent norms.  We then apply Corollary~\ref{cor-generalRKHS-range-space} to  $\H_2=A^2$ and $\H_1=\H(b)=\MM(T_a^*)$, which gives the result.
%

\begin{acknowledgements} We would like to warmly thank the anonymous referee for his/her remarks leading to a real improvement of the paper. In an earlier version, we had stated Theorem~\ref {Thm:invariant-subspace-backward} with the additional hypothesis that for every $\mathcal E\in \mbox{Lat}(X_{\H_1})$, the space $M_{\varphi,\H_2}^*\mathcal E$ is dense in $\mathcal E$. It was the referee's suggestion to use the Sz.-Nagy--Foias theory to obtain that particular property as a consequence of the other hypothesis in the theorem (see Lemma~\ref{lem-densite}). 
\end{acknowledgements}

%
%

\bibliographystyle{mscplain.bst}

\end{document}